\documentclass[11pt,twoside,reqno]{amsart}
\pdfoutput=1
\usepackage{amssymb,amsthm,amsthm,amsfonts,latexsym}
\usepackage{amsmath}
\usepackage{mathrsfs}
\usepackage{stmaryrd}
\usepackage{accents}
\usepackage[latin1]{inputenc}
\usepackage{color}
\usepackage{xspace}
\usepackage{tikz-cd}
\usepackage{float}
\usepackage[]{hyperref}
\usepackage[]{cleveref}

\allowdisplaybreaks

\theoremstyle{plain}
\newtheorem{teorema}{Theorem}[section]
\newtheorem*{teorema*}{Theorem}
\newtheorem{lema}[teorema]{Lemma}
\newtheorem{proposicion}[teorema]{Proposition}
\newtheorem{corolario}[teorema]{Corollary}
\newtheorem*{corolario*}{Corollary}

\theoremstyle{definition}
\newtheorem{definicion}[teorema]{Definition}

\theoremstyle{remark}
\newtheorem{observacion}[teorema]{Remark}

\crefname{teorema}{Theorem}{Theorems}
\crefname{lema}{Lemma}{Lemmas}
\crefname{proposicion}{Proposition}{Propositions}
\crefname{corolario}{Corollary}{Corollaries}
\crefname{definicion}{Definition}{Definitions}
\crefname{observacion}{Remark}{Remarks}
\crefname{pregunta}{Conjecture}{Conjectures}

\newcommand{\st}{\ensuremath{{\mkern2mu\accentset{\circ}{\mathrm{st}}\mkern2mu}}}
\newcommand{\lk}{\ensuremath{\mathrm{lk}}}
\def\co{\colon}
\def\wt{\widetilde}
\newcommand{\tq}{\mathrel{{\ensuremath{\: : \: }}}}
\def\Z{\mathbb{Z}}

\def\Aut{\mathrm{Aut}}
\def\SAut{\mathrm{SAut}}
\def\Cay{\mathrm{Cay}}
\def\PB{\mathrm{PB}}
\def\FC{\mathrm{FC}}
\def\op{\mathrm{op}}
\newcommand{\rk}[1]{\mathrm{rk}\left( #1 \right)}
\def\F{\mathbb{F}}

\def\cc{\bullet}

\def\CondA{(iii)\xspace}
\def\CondB{(i)\xspace}
\def\CondC{(ii)\xspace}

\let\emptyset\varnothing

\begin{document}

\title[The complex of partial bases of a free group]{The complex of partial bases of a free group}

\author[I. Sadofschi Costa]{Iv\'an Sadofschi Costa}

\thanks{The author was partially supported by a CONICET doctoral fellowship and grants CONICET PIP 112-201101-00746 and UBACyT 20020160100081BA}

\address{Universidad de Buenos Aires. Facultad de Ciencias Exactas y Naturales. Departamento de Matemática. Buenos Aires, Argentina.}

\address{CONICET-Universidad de Buenos Aires. Instituto de Investigaciones Matemáticas Luis A. Santaló (IMAS). Buenos Aires, Argentina.}

\email{isadofschi@dm.uba.ar}

\begin{abstract}
We prove that the simplicial complex whose simplices are the nonempty partial bases of $\F_n$ is homotopy equivalent to a wedge of $(n-1)$-spheres. Moreover, we show that it is Cohen-Macaulay.
\end{abstract}

  % Topological methods in group theory}


\maketitle

\section{Introduction}

A \textit{partial basis} of a free group $\F$ is a subset of a basis of $\F$. If $\F$ is a free group of finite rank, we denote the simplicial complex with simplices the nonempty partial bases of $\F$ by $\PB(\F)$. Our main result is \cref{PBFnCohenMacaulay}, which states that $\PB(\F_n)$ is Cohen-Macaulay.
To prove this, we build on previous results on related objects such as the poset $\FC(\F_n)$ of proper free factors of $\F_n$, which was studied by Hatcher and Vogtmann. In \cite{HV}, they proved that its order complex $\mathcal{K}(\FC(\F_n))$ is Cohen-Macaulay. Another closely related complex, defined by Day and Putman, is the complex $\mathcal{B}(\F_n)$ whose simplices are sets $\{C_0,\ldots,C_k\}$ of conjugacy classes of $\F_n$ such that there exists a partial basis $\{v_0,\ldots,v_k\}$ with $C_i=\llbracket v_i \rrbracket $ for $0\leq i \leq k$. Day and Putman prove that $\mathcal{B}(\F_n)$ is $0$-connected for $n\geq 2$ and $1$-connected for $n\geq 3$ \cite[Theorem A]{DP}, that a certain quotient is $(n-2)$-connected \cite[Theorem B]{DP} and they conjecture that $\mathcal{B}(\F_n)$ is $(n-2)$-connected \cite[Conjecture 1.1]{DP}.
The complexes $\PB(\F_n)$, $\mathcal{B}(\F_n)$ and $\FC(\F_n)$ may be seen as analogues for $\Aut(\F_n)$ of the curve complex of a surface. Thus, these objects should give information about $\Aut(\F_n)$. In \cite{DP}, for example, $\mathcal{B}(\F_n)$ is used to prove that the Torelli subgroup is finitely generated.

In Section \ref{SeccionPresentacion}, we obtain a presentation for $\SAut(\F_n,\{v_1,\ldots,v_l\})$, analogous to Gersten's presentation of $\SAut(\F_n)$ which was used by Day and Putman in \cite{DP}. First we use McCool's method to present $\Aut(\F_n,\{v_1,\ldots,v_l\})$ and then we apply the Reidemeister-Schreier method.

In Section \ref{Seccion1Conexion} we prove that the links $\lk(B,\PB(\F_n))$ are $0$-connected if $n-|B|\geq 2$ and $1$-connected if $n-|B|\geq 3$. The proof is the proof of \cite[Theorem A]{DP} with only minor modifications. Instead of Gersten's presentation, we use the presentation obtained in Section \ref{SeccionPresentacion}.

In Section \ref{SeccionQuillen} we prove \cref{QuillenExtendido}, which is a version of a result due to Quillen \cite[Theorem 9.1]{Q}. Our version produces an explicit basis of the top homology group of $X$. The proof is based on Piterman's proof of Quillen's Theorem \cite[Teorema 2.1.28]{P}, which uses Barmak-Minian's non-Hausdorff mapping cylinder argument \cite{BM}.

In Section \ref{SeccionPrincipal} we prove \cref{PBFnCohenMacaulay}. The idea is to compare the link $\lk(B,\PB(\F_n))$ (which is $(n-|B|-1)$-dimensional)  with $\FC(\F_n)_{>\langle B\rangle }$ (which is $(n-|B|-2)$-dimensional). In order do this, we have to consider the $(n-|B|-2)$-skeleton of $\lk(B,\PB(\F_n))$. Finally, using the basis given by \cref{QuillenExtendido}, we can understand what happens when we pass from $\lk(B,\PB(\F_n)^{(n-2)})$ to $\lk(B,\PB(\F_n))$. We proceed by induction on $n-|B|$ and to start we need the result proved in Section \ref{Seccion1Conexion}.

\section[A presentation for SAut(Fn, \{v1, ... , vl\}]{A presentation for $\SAut(\F_n, \{v_1,\ldots, v_l\})$}\label{SeccionPresentacion}

The main result of this section is \cref{teoPresentacion}, which gives a finite presentation for the group $\SAut(\F_n, \{v_1,\ldots, v_l\})$. When $l=0$, this presentation reduces to the presentation of $\SAut(\F_n)$ given by Gersten in \cite{Gersten} and used by Day and Putman in \cite{DP}. To obtain this presentation we first obtain a presentation for $\Aut(\F_n, \{v_1,\ldots, v_l\})$ using McCool's method \cite{MC2} and then we get to the desired presentation using the Reidemeister-Schreier method.

\subsection{Definitions and Notations}
Throughout the paper automorphisms act on the left and compose from right to left as usual.
Let $\F_n$ be a free group with basis $\{v_1,\ldots,v_n\}$. Recall that $\SAut(\F_n)$, \textit{the special automorphism group of $\F_n$}, is the subgroup of $\Aut(\F_n)$ consisting of automorphisms whose images in $\Aut(\Z^n)$ have determinant $1$.
If $A$ is a subset of $\F_n$ we define $\Aut(\F_n, A )= \{ \phi \in\Aut(\F_n) \tq \phi|_A=1_A\}$ and $\SAut(\F_n, A)= \{ \phi \in\SAut(\F_n) \tq \phi|_A=1_A\}$.

Let $L=\{v_1,\ldots,v_n,v_1^{-1},\ldots, v_n^{-1}\}\subseteq \F_n$ be the set of letters. We fix a number $l$, $1\leq l\leq n$ and we define $L'=L-\{v_1,\ldots,v_l,v_1^{-1},\ldots, v_l^{-1}\}$.
We consider the subgroup $\Omega(\F_n)$ of $\Aut(\F_n)$ given by the automorphisms that permute $L$. The order of $\Omega(\F_n)$ is $2^n n!$. 
If $A\subseteq L$, $a\in A$ and $a^{-1}\notin A$, there is an automorphism $(A;a)$ of $\F_n$ defined on $L$ by
$$
(A;a)(x)=\begin{cases}
	      x & \text{ if } x\in\{a,a^{-1}\} \\
	     a^{-1}xa &  \text{ if }x,x^{-1}\in A \\
	     xa  &\text{ if } x\in A, x^{-1}\in L-A\\
	     a^{-1}x & \text{ if  } x^{-1}\in A, x\in L-A \\
	     x &\text{ if } x,x^{-1} \in L- A 	     
         \end{cases}
$$ The set of these automorphisms will be denoted by $\Lambda(\F_n)$.
We consider the set of Whitehead automorphims $\mathcal{W}(\F_n ) = \Lambda(\F_n ) \cup \Omega(\F_n )$.

If $a,b\in L$ with $a\neq b^{\pm 1}$ we denote by $E_{a,b}$ the automorphism that maps $a$ to $ab$ and fixes $L-\{a,a^{-1}\}$ and by $M_{a,b}$ the automorphism that maps $a$ to $ba$ and fixes $L-\{a,a^{-1}\}$. 
We have  $E_{a,b}=(\{a,b\};b)$, $M_{a,b}=(\{a^{-1},b^{-1}\};b^{-1})$ and $E_{a,b}=M_{a^{-1},b^{-1}}$.
If $a,b\in L$ and $a\neq b^{\pm 1}$ we also consider the automorphism $w_{a,b}$ that takes $a$ to $b^{-1}$, $b$ to $a$ and fixes $L-\{a,a^{-1},b,b^{-1}\}$.

\subsection{McCool's method}

We recall the classical presentation of $\Aut(\F_n)$ obtained by McCool.

\begin{teorema}[{\cite{MC1}}]\label{PresentacionMcCool}
There is a presentation of $\Aut(\F_n)$ with generators $\mathcal{W}(\F_n)$ and relators R1-R7 below.

\begin{align}
 (A;a)^{-1}&=(A-\{a\}+\{a\}^{-1};a^{-1}) \tag{R1} \\
     (A;a)(B;a)&=(A\cup B;a) \text{\qquad if $A\cap B= \{a\}$}   \tag{R2} \\
   [(A;a),(B;b)] &=1 \text{\qquad if $A\cap B=\emptyset$, $a^{-1}\notin B$, $b^{-1}\notin A$} \tag{R3} \\
   (B;b)(A;a) &=(A\cup B-\{b\};a)(B;b) \text{\qquad if $A\cap B=\emptyset$, $a^{-1}\notin B$, $b^{-1}\in A$}  \tag{R4} \\
   ( A-\{a\}\cup \{a\}^{-1};b)(A;a) &=(A-\{b\}\cup \{b^{-1}\};a)w_{a,b} \text{\qquad if $b\in A$, $b^{-1}\notin A$, $a\neq b$}  \tag{R5} \\
  T(A;a)T^{-1}&=(T(A);T(a)) \text{\qquad if $T\in \Omega(\F_n)$} \tag{R6} \\
     &\text{Multiplication table of $\Omega(\F_n)$}  \tag{R7}
\end{align}
Additionally the following relations hold
\begin{align}
    (A;a) =(L-A; a^{-1} ) (L-\{a^{-1}\};a) &=   (L-\{a^{-1}\}; a) (L-A; a^{-1})            \tag{R8} \\
    (L-\{b^{-1}\};b)(A;a)(L-\{b\};b^{-1})&=(A;a)                 \text{\qquad if $b,b^{-1}\in L-A$} \tag{R9} \\
     (L-\{b^{-1}\};b)(A;a)(L-\{b\};b^{-1})&=(L-A;a^{-1})          \text{\qquad if $b\neq a$, $b\in A$, $b^{-1}\in L-A$}\tag{R10} 
\end{align}

\end{teorema}

The \textit{length} of an element $u\in \F_n$, denoted $|u|$, is the number of letters of the unique reduced word in $L$ that represents $u$. The \textit{total length} of an $m$-tuple $(u_1,\ldots,u_m)$ of elements of $\F_n$ is $|u_1|+\ldots+|u_m|$. 

\begin{teorema}[{McCool's method, \cite[Section 4. (1)-(2)]{MC2}}]\label{McCoolsMethod}
Let $\F_n$ be the free group with basis $\{v_1,\ldots,v_n\}$. The group $\mathcal{A}=\Aut(\F_n)$ acts on $(\F_n)^{m}$ by $\phi\cdot (u_1,\ldots, u_m)=(\phi(u_1),\ldots,\phi(u_m))$.
Then if $U=(u_1,\ldots, u_m)$ is an $m$-tuple of words in $\F_n$, the stabilizer $\mathcal{A}_U$ of $U$ is finitely presented.

Moreover, we can construct a finite $2$-complex $K$ with fundamental group $\mathcal{A}_U$ as follows. Let $K^{(0)}$ be the set of tuples in the orbit of $U$ which have minimum total length. We may assume $U\in K^{(0)}$. Each triple $(V,V', \phi )$ such that $V,V'\in K^{(0)}$, $\phi\in\mathcal{W}(\F_n)$ and $\phi(V)=V'$ represents a directed edge from $V$ to $V'$ labelled $\phi$. The directed edges $(V,V', \phi )$ and $(V',V,\phi^{-1})$ together determine a single $1$-cell in $K^{(1)}$. Finally a $2$-cell is attached following each closed edge path in $K^{(1)}$ such that the word obtained reading the labels is a relator corresponding to a relation of type R1-R10. Then $\mathcal{A}_U=\pi_1(K,U)$.
\end{teorema}

\begin{observacion}
Let $\mathcal{P}$ be the presentation in \cref{PresentacionMcCool} and let $K_\mathcal{P}$ be the associated $2$-complex. The complex $K$ in \cref{McCoolsMethod} is a subcomplex of the covering space of $K_\mathcal{P}$ that corresponds to the subgroup $\mathcal{A}_U\leq \pi_1(K_\mathcal{P},x_0)$.
\end{observacion}

\begin{teorema}\label{Presentacion1}
There is a presentation of $\Aut(\F_n,\{v_1,\ldots,v_l\})$ with generators $$\mathcal{W}(\F_n) \cap \Aut(\F_n,\{v_1,\ldots,v_l\})$$ and relators R1-R7 that involve only those generators.
\begin{proof}
Let $U=(v_1,\ldots,v_l)$. We use \cref{McCoolsMethod} to construct a $2$-complex $K$ with fundamental group $\Aut(\F_n,\{v_1,\ldots,v_l\})=\mathcal{A}_U$. We note that $U$ has minimum total length and the $0$-skeleton $K^{(0)}$ is the set of tuples $V=(v_{\sigma(1)}^{s_1}, \ldots, v_{\sigma(l)}^{s_l})$ with $\sigma \in S_n$ and $s_i\in\{1,-1\}$.

Now we can obtain a presentation of $\pi_1(K,U)$. The presentation has a generator for each edge and a relation for each $2$-cell of $K$ and also a relation for each edge in a fixed spanning tree of $K^1$. The spanning tree we choose consists of an edge $(V,U,\sigma_V)$ for each $V\in K^0-\{U\}$, for some fixed $\sigma_V\in \Omega(\F_n)$.
We note that if $x\in L$ and $(A;a)(x)\in L$ then $x=(A;a)(x)$. Therefore an edge labeled $(A;a)$ is necessarily a loop. Hence the relations of types R1-R5 and R8-R10 are products of loops on a same vertex $V$. If we consider such a relation
$$(V,V,(A_1;a_1))\cdots (V,V,(A_k;a_k))$$
with $V\neq U$, using relations R6 we can replace it by
$$(V,U,\sigma)^{-1}(U,U,(\sigma(A_1);\sigma(a_1)))\cdots (U,U,(\sigma(A_k);\sigma(a_k)))(V,U,\sigma)$$
for some $\sigma \in \Omega(\F_n)$, which in turn is equivalent to
$$(U,U,(\sigma(A_1);\sigma(a_1)))\cdots (U,U,(\sigma(A_k);\sigma(a_k)))$$
But this relation already appears in the presentation. For this reason we can discard every relation of types R1-R5 and R8-R10 which is based at a vertex different from $U$. 
Now the generators $(V,V,(A;a))$ with $V\neq U$ appear only in the relations of type R6. We will show that we can remove almost all relations of type R6, so that the only ones left are those based at $U$ and those given by
$$(V,V,(A;a))=(V,U,\sigma_V)^{-1} (U,U, (\sigma_V(A);\sigma_V(a))) (V,U,\sigma_V)$$
with $V\neq U$. Then we will eliminate the generators $(V,V,(A;a))$ along with these relations, so that the only relations of type R6 left are based at $U$.

First we show that relations 
$$(V,V,(A;a))=(V,V',\tau)^{-1}(V',V',(\tau(A);\tau(a)))(V,V',\tau)$$
with $V,V'\neq U$ and $\tau(V)=V'$ are redundant. To do this we take $\sigma\in \Omega(\F_n)$ such that $\sigma(U)=V$. Conjugating by $(U,V,\sigma)$ and using a relation of type R7 we can replace our relation by 
$$(U,V,\sigma)^{-1}(V,V,(A;a))(U,V,\sigma)=(U,V',\tau\sigma)^{-1}(V',V',(\tau(A);\tau(a)))(U,V',\tau\sigma)$$
Now using $(U,U,(\sigma^{-1}(A);\sigma^{-1}(a)))=(U,V,\sigma)^{-1}(V,V,(A;a))(U,V,\sigma)$ we see that this is equivalent to 
$$(U,U,(\sigma^{-1}(A);\sigma^{-1}(a)))=(U,V',\tau\sigma)^{-1}(V',V',(\tau(A);\tau(a)))(U,V',\tau\sigma)$$
which is repeated.

Now we show that if $\tau \neq \sigma_V$ satisfies $\tau(V)=U$, the relation
$$(V,V,(A;a))=(V,U,\tau)^{-1} (U,U, (\tau(A);\tau(a))) (V,U,\tau)$$
is redundant, since it can be replaced by
$$(V,U,\sigma_V)^{-1} (U,U, (\sigma_V(A);\sigma_V(a))) (V,U,\sigma_V) = (V,U,\tau)^{-1} (U,U, (\tau(A);\tau(a))) (V,U,\tau)$$
that can be rewritten (using R7) as
$$(U,U, (\sigma_V(A);\sigma_V(a)))  = (U,U,\tau\sigma_V^{-1})^{-1} (U,U, (\tau(A);\tau(a))) (U,U,\tau\sigma_V^{-1})$$
that is a relation of type R6 based at $U$. Relations $$(U,U,(A;a))=(U,V,\tau)^{-1} (V,V, (\tau(A);\tau(a))) (U,V,\tau)$$ are obviously equivalent to those considered before.

In this way we can eliminate the generators corresponding to the edges $(V,V,(A;a))$ and the only relations of type R6 left are those given by loops in $U$.
At this point the only generators left are those labelled with an element of  $\Omega(\F_n)$ or based at $U$ and the only relators left are those of type R7 or of types R1-R10 based at $U$.
In a similar way we eliminate the generators that are not based at $U$  along with the relations of type R7 not based at $U$.

To finish we must eliminate the relators of type R8-R10 using that they follow from R1-R7 (see \cite[3.]{MC1}). If $l=0$ we already know that these relations are redundant. If $l>1$, the only letters fixed by $(L-\{a^{-1}\};a)$ are $a$ and $a^{-1}$ so in this case there are no relators of type R8-R10. To deal with the case $l=1$ we must check that the every relation of type R1-R7 used in \cite[3.]{MC1} to check the relation of type R8-R10 we intend to eliminate fixes $v_1$. For example, if we are eliminating a relation of type R9, we have $b=v_1^{\pm 1}$ so any generator $(X;b)$ or $(X;b^{-1})$ fixes $v_1$. In addition $(A;a)$ fixes $v_1$. So every generator in the intermediate steps fixes $v_1$ and we are done.
\end{proof}
\end{teorema}

\subsection{The Reidemeister-Schreier method}
To simplify the presentation we change the generating set following Gersten \cite{Gersten}.

\begin{teorema}\label{PresentacionAutFijando}
The group $\Aut(\F_n,\{v_1,\ldots,v_l\})$ has a presentation with generators $$\{ M_{a,b}\tq  a\in L', b\in L, a\neq b^{\pm 1} \}\cup (\Omega(\F_n)\cap \Aut(\F_n,\{v_1,\ldots,v_l\}))$$ subject to the following relations:
\begin{itemize} 
 \item[S0.] Multiplication table of $\Omega(\F_n)\cap \Aut(\F_n,\{v_1,\ldots,v_l\})$.
 
 \item[S1.] $M_{a,b}M_{a,b^{-1}}=1$.
 
 \item[S2.]  $[M_{a,b},M_{c,d}]=1$ if $b\neq c^{\pm 1}$, $a\neq d^{\pm 1}$ and $a\neq c$.
 
 \item[S3.] $[M_{b,a^{-1}},M_{c,b^{-1}}]=M_{c,a}$.
 
 \item[S4.] $w_{a,b}= M_{b^{-1},a^{-1}}M_{a^{-1},b} M_{b,a}$.
 
 \item[S5.] $\sigma M_{a,b} \sigma^{-1} = M_{\sigma(a),\sigma(b)}$ if $\sigma \in \Omega(\F_n)\cap \Aut(\F_n,\{v_1,\ldots,v_l\})$.
 \end{itemize}
 
 \begin{observacion}
 In S5 both the generator and the automorphism are denoted by $\sigma$. We will do this repeatedly.
 \end{observacion}

 \begin{proof}
 The generators of this presentation are elements of $\Aut(\F_n,\{v_1,\ldots,v_l\})$ that verify relations S0-S5. We have
 $$(A;a)= \prod_{b\in A,b\neq a} E_{b,a} = \prod_{b\in A,b\neq a} M_{b^{-1},a^{-1}}$$ (S2 guarantees the product is well-defined). Therefore, by \cref{Presentacion1} we have a generating set of $\Aut(\F_n,\{v_1,\ldots,v_l\})$.
 It suffices to check that R1-R7 can be deduced from S0-S5. The proof is the same as in  \cite[Theorem 1.2]{Gersten}.
 \end{proof}
\end{teorema}

\begin{teorema}[{c.f. \cite[Theorem 1.4]{Gersten}}]\label{PresentacionRS}
The group $\SAut(\F_n,\{v_1,\ldots,v_l\})$ has a presentation with generators $$\{ M_{a,b}\tq  a\in L', b\in L,  a\neq b^{\pm 1}  \}\cup (\Omega(\F_n)\cap \SAut(\F_n,\{v_1,\ldots,v_l\}))$$ subject to the following relations:
 \begin{itemize}
 \item[S0.] Multiplication table of $\Omega(\F_n)\cap \SAut(\F_n,\{v_1,\ldots,v_l\})$.
 
 \item[S1.] $M_{a,b}M_{a,b^{-1}}=1$.
 
 \item[S2.]  $[M_{a,b},M_{c,d}]=1$ if $b\neq c^{\pm 1}$, $a\neq d^{\pm 1}$ and $a\neq c$.
 
 \item[S3.] $[M_{b,a^{-1}},M_{c,b^{-1}}]=M_{c,a}$.
 
 \item[S4.] $w_{a,b}= M_{b^{-1},a^{-1}}M_{a^{-1},b} M_{b,a}$.
 
 \item[S5.] $\sigma M_{a,b} \sigma^{-1} = M_{\sigma(a),\sigma(b)}$ if $\sigma \in \Omega(\F_n)\cap \SAut(\F_n,\{v_1,\ldots,v_l\})$.
\end{itemize}

 \begin{proof}
 The presentation is obtained applying the Reidemeister-Schreier method to the presentation $\mathcal{P}$ of \cref{PresentacionAutFijando}. That is, we consider the associated $2$-complex $K_\mathcal{P}$ and we construct the covering space corresponding to the subgroup $\SAut(\F_n,\{v_1,\ldots,v_l\})\leq \Aut(\F_n,\{v_1,\ldots,v_l\})$ which has index $2$. This covering has two cells over each cell of $K_\mathcal{P}$. We obtain a presentation of its fundamental group and then we eliminate generators and relators until we get the desired presentation.
 \end{proof}
\end{teorema}

\begin{lema}[{\cite[Lemma 1.3]{Gersten} }]
The group $\Omega(\F_n)\cap \SAut(\F_n)$ has a presentation with generators $\{ w_{a,b}\tq a,b\in L, a\neq b^{\pm 1} \}$ and relations:
\begin{itemize}
 \item $w_{a,b^{-1}}=w_{a,b}^{-1}$
 \item $w_{a,b}w_{c,d} w_{a,b}^{-1} = w_{w_{a,b}(c),w_{a,b}(d)}$
 \item $w_{a,b}^4=1$
\end{itemize}
\end{lema}

\begin{corolario}\label{PresentacionOmegaCapSAut}
The group $\Omega(\F_n)\cap \SAut(\F_n,\{v_1,\ldots,v_l\})$ has a presentation with generators 
$\{ w_{a,b}\tq a,b\in L', a\neq b^{\pm 1} \}$ and the following relations:
\begin{itemize}
 \item $w_{a,b^{-1}}=w_{a,b}^{-1}$
 \item $w_{a,b}w_{c,d} w_{a,b}^{-1} = w_{w_{a,b}(c),w_{a,b}(d)}$
 \item $w_{a,b}^4=1$
\end{itemize}
\begin{proof}
This follows immediately from the previous lemma using $\Omega(\F_n)\cap \SAut(\F_n,\{v_1,\ldots,v_l\})\simeq \Omega(\F_{n-l})\cap \SAut(\F_{n-l})$.
\end{proof}
\end{corolario}

\begin{teorema}\label{teoPresentacionPrevio}
 The group $\SAut(\F_n,\{v_1,\ldots,v_l\})$ has a presentation with generators
  $$\{ M_{a,b}\tq a\in L', b\in L,  a\neq b^{\pm 1}  \}\cup \{ w_{a,b}\tq a,b\in L',  a\neq b^{\pm 1}  \}$$ subject to the following relations
 	\begin{itemize}
	  \item[(1)] $M_{a,b}M_{a,b^{-1}}=1$.
	  \item[(2)]  $[M_{a,b},M_{c,d}]=1$ if $b\neq c^{\pm 1}$, $a\neq d^{\pm 1}$ and $a\neq c$.
	  \item[(3)] $[M_{b,a^{-1}},M_{c,b^{-1}}]=M_{c,a}$.
	  \item[(4)] $w_{a,b}=M_{b^{-1},a^{-1}}M_{a^{-1},b}M_{b,a}$.
	  \item[(5')$\!$] $w_{a,b}M_{c,d} w_{a,b}^{-1} = M_{w_{a,b}(c),w_{a,b}(d)}$.
	  \item[(6)] $w_{a,b}^4=1$.
	\end{itemize}
 \begin{proof}
This presentation can be obtained from \cref{PresentacionRS} using \cref{PresentacionOmegaCapSAut}.
 \end{proof}
\end{teorema}

Now we state and prove the main result of this section.

\begin{teorema}\label{teoPresentacion}
If $n-l\geq 3$, the group $\SAut(\F_n,\{v_1,\ldots,v_l\})$ has a presentation with generators 
$$\{ M_{a,b}\tq a\in L', b\in L, a\neq b^{\pm 1}  \}\cup \{ w_{a,b}\tq a,b\in L',   a\neq b^{\pm 1}  \}$$ subject to the following relations
subject to the following relations:
 \begin{enumerate}
	  \item $M_{a,b}M_{a,b^{-1}}=1$.
	  \item  $[M_{a,b},M_{c,d}]=1$ if $b\neq c^{\pm 1}$, $a\neq d^{\pm 1}$ and $a\neq c$.
	  \item $[M_{b,a^{-1}},M_{c,b^{-1}}]=M_{c,a}$.
	  \item $w_{a,b}=M_{b^{-1},a^{-1}}M_{a^{-1},b}M_{b,a}$.
  	  \item $w_{a,b}=w_{a^{-1},b^{-1}}$.
	  \item $w_{a,b}^4=1$.
\end{enumerate}
 \begin{proof}
 We only need to show that relation (5') follows from relations (1)-(6). To do this we do the same computations as in \cite[Proof of Theorem 2.7]{Gersten}.
We separate in cases. In each case different letters represent elements which belong to different orbits of the action of $\Z_2$ on $L$ given by $x\mapsto x^{-1}$.

\begin{itemize}
\item Cases with four orbits.
\begin{itemize}
\item  $w_{a,b}M_{c,d} w_{a,b}^{-1} = M_{w_{a,b}(c),w_{a,b}(d)}$

\begin{align*}
 w_{a,b}M_{c,d} w_{a,b}^{-1} &= M_{b^{-1},a^{-1}}M_{a^{-1},b}M_{b,a} M_{c,d} w_{a,b}^{-1} \\
\text{(Using 2)} &= M_{b^{-1},a^{-1}}M_{a^{-1},b} M_{c,d}M_{b,a}  w_{a,b}^{-1} \\
\text{(Using 2)} &= M_{b^{-1},a^{-1}}M_{c,d}M_{a^{-1},b} M_{b,a}  w_{a,b}^{-1} \\
\text{(Using 2)} &= M_{c,d} M_{b^{-1},a^{-1}}M_{a^{-1},b} M_{b,a} w_{a,b}^{-1} \\
&= M_{c,d}w_{a,b} w_{a,b}^{-1} \\
&= M_{c,d}\\
&= M_{w_{a,b}(c),w_{a,b}(d)}
\end{align*}

\end{itemize}
\item Cases with three orbits.

In every case we use $w_{a,b}=w_{a^{-1},b^{-1}}=M_{b,a}M_{a,b^{-1}}M_{b^{-1},a^{-1}}$.

\begin{itemize}
%\item[] $\langle a\rangle =\langle c\rangle$

\item  $w_{a,b}M_{a,d} w_{a,b}^{-1} = M_{w_{a,b}(a),w_{a,b}(d)}$

\begin{align*}
w_{a,b}M_{a,d} w_{a,b}^{-1} &= M_{b,a}M_{a,b^{-1}}M_{b^{-1},a^{-1}} M_{a,d} M_{b^{-1},a} M_{a,b} M_{b,a^{-1}} \\
\text{(Using 3)}  &= M_{b,a}M_{a,b^{-1}} M_{b^{-1},d}M_{a,d}M_{b^{-1},a^{-1}}  M_{b^{-1},a} M_{a,b} M_{b,a^{-1}} \\
\text{(Using 1)}  &= M_{b,a}M_{a,b^{-1}} M_{b^{-1},d}M_{a,d} M_{a,b} M_{b,a^{-1}} \\
\text{(Using 2)}  &= M_{b,a}M_{a,b^{-1}} M_{a,d} M_{b^{-1},d} M_{a,b} M_{b,a^{-1}} \\
\text{(Using 3)}  &= M_{b,a}M_{a,b^{-1}} M_{a,d} M_{a,d^{-1}} M_{a,b} M_{b^{-1},d} M_{b,a^{-1}} \\
\text{(Using 1)}  &= M_{b,a}M_{a,b^{-1}}  M_{a,b} M_{b^{-1},d} M_{b,a^{-1}} \\
\text{(Using 1)}  &= M_{b,a} M_{b^{-1},d} M_{b,a^{-1}} \\
\text{(Using 2)}  &= M_{b,a} M_{b^{-1},d} M_{b,a^{-1}} \\
\text{(Using 2)}  &=  M_{b^{-1},d}M_{b,a} M_{b,a^{-1}} \\
\text{(Using 1)}  &=  M_{b^{-1},d} \\
&= M_{w_{a,b}(a),w_{a,b}(d) } \\
\end{align*}

\item  $w_{a,b}M_{a^{-1},d} w_{a,b}^{-1} = M_{w_{a,b}(a^{-1}),w_{a,b}(d)}$
 
 By (5) this case follows from the previous one.

%\item[] $<a>=<d>$

\item  $w_{a,b}M_{c,a} w_{a,b}^{-1} = M_{w_{a,b}(c),w_{a,b}(a)}$

By (5) this case follows from the next one.

\item  $w_{a,b}M_{c,a^{-1}} w_{a,b}^{-1} = M_{w_{a,b}(c),w_{a,b}(a^{-1})}$

\begin{align*}
w_{a,b}M_{c,a^{-1}} w_{a,b}^{-1} &= M_{b,a}M_{a,b^{-1}}M_{b^{-1},a^{-1}} M_{c,a^{-1}} M_{b^{-1},a} M_{a,b} M_{b,a^{-1}} \\
\text{(Using 2)} &=  M_{b,a}M_{a,b^{-1}}M_{b^{-1},a^{-1}}  M_{b^{-1},a}  M_{c,a^{-1}} M_{a,b} M_{b,a^{-1}} \\
\text{(Using 1)} &=  M_{b,a}M_{a,b^{-1}}  M_{c,a^{-1}} M_{a,b} M_{b,a^{-1}} \\
\text{(Using 3)} &=  M_{b,a}M_{c,b}M_{c,a^{-1}}M_{a,b^{-1}} M_{a,b} M_{b,a^{-1}} \\
\text{(Using 1)} &=  M_{b,a}M_{c,b}M_{c,a^{-1}} M_{b,a^{-1}} \\
\text{(Using 2)} &=  M_{b,a}M_{c,b} M_{b,a^{-1}} M_{c,a^{-1}} \\
\text{(Using 3)} &=  M_{c,b} M_{b,a} M_{c,a} M_{b,a^{-1}} M_{c,a^{-1}} \\
\text{(Using 2)} &= M_{c,b} \\
&= M_{w_{a,b}(c),w_{a,b}(a^{-1}) } \\
 \end{align*}

%\item[] $<b>=<c>$

\item  $w_{a,b}M_{b,d} w_{a,b}^{-1} = M_{w_{a,b}(b),w_{a,b}(d)}$

\begin{align*}
w_{a,b}M_{b,d} w_{a,b}^{-1} &= M_{b,a}M_{a,b^{-1}}M_{b^{-1},a^{-1}} M_{b,d} M_{b^{-1},a} M_{a,b} M_{b,a^{-1}} \\
\text{(Using 2)} &=  M_{b,a}M_{a,b^{-1}}M_{b^{-1},a^{-1}}  M_{b^{-1},a}  M_{b,d} M_{a,b} M_{b,a^{-1}} \\
\text{(Using 1)} &=  M_{b,a}M_{a,b^{-1}}  M_{b,d} M_{a,b} M_{b,a^{-1}} \\
\text{(Using 3)} &=  M_{b,a}M_{a,b^{-1}}  M_{a,b} M_{b,d} M_{a,d} M_{b,a^{-1}} \\
\text{(Using 1)} &=  M_{b,a} M_{b,d} M_{a,d} M_{b,a^{-1}} \\
\text{(Using 3)} &=  M_{b,a} M_{b,d} M_{b,d^{-1}} M_{b,a^{-1}} M_{a,d} \\
\text{(Using 1)} &=  M_{b,a}  M_{b,a^{-1}} M_{a,d} \\
\text{(Using 1)} &= M_{a,d} \\
&= M_{w_{a,b}(a),w_{a,b}(d) } \\
\end{align*}

\item  $w_{a,b}M_{b^{-1},d} w_{a,b}^{-1} = M_{w_{a,b}(b^{-1}),w_{a,b}(d)}$

 By (5) this case follows from the previous one.

%\item[] $<b>=<d>$

\item  $w_{a,b}M_{c,b} w_{a,b}^{-1} = M_{w_{a,b}(c),w_{a,b}(b)}$

\begin{align*}
 w_{a,b}M_{c,b} w_{a,b}^{-1} &= M_{b,a}M_{a,b^{-1}}M_{b^{-1},a^{-1}} M_{c,b} M_{b^{-1},a} M_{a,b} M_{b,a^{-1}} \\
\text{(Using 3)}  &= M_{b,a}M_{a,b^{-1}} M_{c,a} M_{c,b} M_{b^{-1}, a^{-1}}  M_{b^{-1},a} M_{a,b} M_{b,a^{-1}} \\
\text{(Using 1)}  &= M_{b,a}M_{a,b^{-1}} M_{c,a} M_{c,b} M_{a,b} M_{b,a^{-1}} \\
\text{(Using 2)}  &= M_{b,a}M_{a,b^{-1}} M_{c,a} M_{a,b}  M_{c,b}M_{b,a^{-1}} \\
\text{(Using 3)}  &= M_{b,a}M_{c,a} M_{a,b^{-1}}M_{c,b^{-1}}  M_{a,b}  M_{c,b}M_{b,a^{-1}} \\
\text{(Using 2)}  &= M_{b,a}M_{c,a} M_{a,b^{-1}}  M_{a,b}M_{c,b^{-1}}  M_{c,b}M_{b,a^{-1}} \\
\text{(Using 1)}  &= M_{b,a}M_{c,a} M_{c,b^{-1}}  M_{c,b}M_{b,a^{-1}} \\
\text{(Using 1)}  &= M_{b,a}M_{c,a} M_{b,a^{-1}} \\
\text{(Using 2)}  &= M_{c,a}M_{b,a} M_{b,a^{-1}} \\
\text{(Using 1)}  &= M_{c,a} \\
&= M_{w_{a,b}(c),w_{a,b}(b) } \\
 \end{align*}

\item  $w_{a,b}M_{c,b^{-1}} w_{a,b}^{-1} = M_{w_{a,b}(c),w_{a,b}(b^{-1})}$

 By (5) this case follows from the previous one.

\end{itemize}

\item Cases with two orbits

\begin{itemize}
%\item[] $<a>=<c>$ y $<b>=<d>$

\item  $w_{a,b}M_{a,b} w_{a,b}^{-1} = M_{w_{a,b}(a),w_{a,b}(b)}$

We consider $c\in L'$ such that the orbit of $c$ is different from the orbits of $a$ and $b$. Using $M_{a,b}=[M_{c,b^{-1}},M_{a,c^{-1}}]$ we have
\begin{align*}
 w_{a,b}M_{a,b}w_{a,b}^{-1}&=  w_{a,b}[M_{c,b^{-1}}, M_{a,c^{-1}}]w_{a,b}^{-1}\\
  &= [w_{a,b}M_{c,b^{-1}}w_{a,b}^{-1}, w_{a,b}M_{a,c^{-1}}w_{a,b}^{-1}] \\
 \text{(Case of 3 orbits already proved)}&= [ M_{w_{a,b}(c),w_{a,b}(b^{-1})}, M_{w_{a,b}(a),w_{a,b}(c^{-1})}] \\
 \text{(Using 3)} &= M_{w_{a,b}(a),w_{a,b}(b)}
\end{align*}

\item  $w_{a,b}M_{a,b^{-1}} w_{a,b}^{-1} = M_{w_{a,b}(a),w_{a,b}(b^{-1})}$

This follows from the previous case taking inverse at both sides.
\item  $w_{a,b}M_{a^{-1},b} w_{a,b}^{-1} = M_{w_{a,b}(a^{-1}),w_{a,b}(b)}$

Using (5) this is equivalent to $w_{a^{-1},b^{-1}} M_{a^{-1},b} w_{a^{-1},b^{-1}}^{-1} = M_{w_{a^{-1},b^{-1}}(a^{-1}),w_{a^{-1},b^{-1}}(b)}$. But this is the previous case applied to $a^{-1}$ and $b^{-1}$.

\item  $w_{a,b}M_{a^{-1},b^{-1} } w_{a,b}^{-1} = M_{w_{a,b}(a^{-1}),w_{a,b}(b^{-1})}$

This follows from the previous case taking inverse at both sides.

%\item[] $<a>=<d>$ y $<b>=<c>$

\item  $w_{a,b}M_{b,a} w_{a,b}^{-1} = M_{w_{a,b}(b),w_{a,b}(a)}$

We consider $c\in L'$ such that the orbit of $c$ is different from the orbits of $a$ and $b$. Using $M_{b,a}=[M_{c,a^{-1}},M_{b,c^{-1}}]$ we have
\begin{align*}
 w_{a,b}M_{b,a}w_{a,b}^{-1}&=  w_{a,b}[M_{c,a^{-1}}, M_{b,c^{-1}}]w_{a,b}^{-1}\\
 &= [w_{a,b}M_{c,a^{-1}}w_{a,b}^{-1}, w_{a,b}M_{b,c^{-1}}w_{a,b}^{-1}] \\
 \text{(Case of 3 orbits already proved)}&= [ M_{w_{a,b}(c),w_{a,b}(a^{-1})}, M_{w_{a,b}(b),w_{a,b}(c^{-1})}] \\
 &= M_{w_{a,b}(b),w_{a,b}(a)}
\end{align*}

\item  $w_{a,b}M_{b,a^{-1}} w_{a,b}^{-1} = M_{w_{a,b}(b),w_{a,b}(a^{-1})}$

This follows from the previous case taking inverse at both sides.

\item  $w_{a,b}M_{b^{-1},a} w_{a,b}^{-1} = M_{w_{a,b}(b^{-1}),w_{a,b}(a)}$

Using (5) this is equivalent to  $w_{a^{-1},b^{-1}}M_{b^{-1},a} w_{a^{-1},b^{-1}}^{-1} = M_{w_{a^{-1},b^{-1}}(b^{-1}),w_{a^{-1},b^{-1}}(a)}$. But this is the previous case applied to $a^{-1}$ and $b^{-1}$.

\item  $w_{a,b}M_{b^{-1},a^{-1}} w_{a,b}^{-1} = M_{w_{a,b}(b^{-1}),w_{a,b}(a^{-1})}$

This follows from the previous case taking inverse at both sides.
\end{itemize}
\end{itemize}

 \end{proof}
\end{teorema}
\section[The links are 1-connected]{The links are $1$-connected}\label{Seccion1Conexion}

Recall that if $G$ is a group and $\mathcal{S}\subseteq G$ is a generating set, the \textit{Cayley graph} $\Cay(G,\mathcal{S})$ is the graph with vertex set $G$ and an edge  $\{g,g'\}$  whenever there exists $s\in\mathcal{S}$ such that $g'=gs$.

\begin{teorema}[c.f. {\cite[Theorem A]{DP}}]\label{linksimplementeconexo}
Let $B$ be a partial basis of $\F_n$. The complex $\lk(B,\PB(\F_n))$ is  connected for $n-|B|\geq 2$ and $1$-connected for $n-|B|\geq 3$.
 \begin{proof}
  The same proof of \cite[Theorem A]{DP} works with subtle changes. We extend $B=\{v_1,\ldots,v_l\}$ to a basis $\{v_1,\ldots, v_n\}$ of $\F_n$. Let $\mathcal{P}=\langle \mathcal{S} \mid \mathcal{R}\rangle$ be the presentation from \cref{teoPresentacion}. We know that $\mathcal{S}$ is a generating set for  $\SAut(\F_n, B)$, even if $n-|B|=2$ (for example, by \cref{teoPresentacionPrevio}).
  
  We define a cellular $\SAut(\F_n,B)$-equivariant map  $\Phi\co \Cay(\SAut(\F_n,B),\mathcal{S})\to \lk(B,\PB(\F_n))$. If $f$ is a $0$-cell, $\Phi(f)=f(v_{l+1})$.
  Now we have to define $\Phi$ on the $1$-cells. If $f - fs$ is a $1$-cell there are three cases:
  \begin{itemize}
   \item $v_{l+1}=s(v_{l+1})$. In this case $\Phi$ maps the entire $1$-cell to $f(v_{l+1})$.
   \item $B\cup \{v_{l+1}, s(v_{l+1})\}$ is a partial basis. In this case $\Phi$ maps the $1$-cell homeomorphically to the edge $\{f(v_{l+1}), fs(v_{l+1}))\}$.
   \item $s=M_{ v_{l+1}^{e},v_i^{e'}  }^{e''}$ for some $1\leq i\leq l$ and $e,e',e''\in\{1,-1\}$. In this case $v_{l+2}=s(v_{l+2})$ and the image of the $1$-cell is the edge-path $f(v_{l+1})-f(v_{l+2})-fs(v_{l+1})$.
  \end{itemize}
  Note that this is well-defined, since the definitions for $f - fs$ and $fs - f$ agree.

  If $n-|B|\geq 2$, the action of $\SAut(\F_n,B)$ on the vertex set of $\lk(B,\PB(\F_n))$ is transitive and the image of $\Phi$ contains every vertex of $\lk(B,\PB(\F_n))$. Since $\Cay(\SAut(\F_n,B),\mathcal{S})$ is connected we conclude that $\lk(B,\PB(\F_n))$ is connected when $n-|B|\geq 2$.

  Now suppose $n-|B|\geq 3$. To prove that $\lk(B,\PB(\F_n))$ is $1$-connected, we will show that $$\Phi_*:\pi_1(\Cay(\SAut(\F_n,B),\mathcal{S}),1)\to \pi_1(\lk(B,\PB(\F_n)),v_{l+1})$$ is surjective and has trivial image.

\textbf{Claim 1.} \textit{The map $\Phi_*\co \pi_1(\Cay(\SAut(\F_n,B),\mathcal{S}),1)\to \pi_1(\lk(B,\PB(\F_n)),v_{l+1})$ has trivial image.}

We will show that $\Phi$ extends to the universal cover $\wt{X}_{\mathcal{P}}$ of $X_\mathcal{P}$. To prove this it is enough to prove that for every relation in $\mathcal{R}$, the image by $\Phi$ of every attaching map (in $\wt{X}_{\mathcal{P}}$) associated to that relation is null-homotopic. Since $\Phi$ is equivariant, it is enough to prove this for the lifts at $1$.

Let $s_1\cdots s_k=1$  be a relation in $\mathcal{R}$. The associated edge-path loop is
$$1-s_1 - s_1 s_2 - \cdots - s_1s_2\cdots s_k=1$$
and its image by $\Phi_*$ is the concatenation of the paths
$$\Phi_*( s_1\cdots s_{i-1} - s_1\cdots s_{i-1}s_i ) = s_1\cdots s_{i-1} \Phi_*( 1- s_i )$$
for $i=1,\ldots, k$.

Inspecting the relations in $\mathcal{R}$ we see that there exists $x\in \{v_{l+1},\ldots, v_n\}$ such that $s_i(x)=x$ for $1\leq i \leq k$ (here we use $n-|B|\geq 3$). Therefore we have that either $v_{l+1}$ equals $x$ or these vertices are joined by an edge. Hence, if $1\leq i \leq k$, $s_1\cdots s_i(v_{l+1})$ and $s_1\cdots s_i(x)=x$ are either equal or joined by an edge. Therefore it suffices to show that the edge-path loop $$x - s_1\cdots s_{i-1} \Phi_*( 1 - s_i ) - x$$ is trivial for every $i$. Using the action of $\SAut(\F_n,B)$, it is enough to show the loops $$x- \Phi_*( 1 - s_i ) - x$$ are trivial. We separate in cases:
\begin{itemize}
 \item If  $v_{l+1}=s_i(v_{l+1})$ it is immediate.
 \item If $B\cup \{v_{l+1}, s(v_{l+1})\}$ is a partial basis we have two cases.
	\begin{itemize}
	  \item If $B\cup \{v_{l+1},s_i(v_{l+1}),x\}$ is a partial basis it is immediate.
	  \item If $B\cup \{v_{l+1},s_i(v_{l+1}),x\}$ is not a partial basis, inspecting $\mathcal{S}$, we see that $s_i=M_{v_{l+1}^{e},x^{e'}}^{e''}$ for certain $e,e',e''\in\{1,-1\}$. Considering $y\in \{v_{l+1},\ldots, v_{n}\}$ distinct from $x$ and $v_{l+1}$ we conclude that the loop $x-v_{l+1} - s_i(v_{l+1}) -x$ contracts to $y$.
	\end{itemize}
 \item If $s_i=M_{ v_{l+1}^{e},v_i^{e'}  }^{e''}$ for some $1\leq i\leq l$ and $e,e',e''\in\{1,-1\}$, we have to show that the loop $x - v_{l+1}- v_{l+2}- s_i(v_{l+1})-x$ is null-homotopic. Again we have two cases.
 	\begin{itemize}
	  \item If $x=v_{l+2}$ it is immediate.
	  \item If $x\neq v_{l+2}$, then $\{x,v_{l+1},v_{l+2} \}$ is a $2$-simplex. Additionally $s_i(v_{l+1})=(v_i^{e'}v_{l+1}^{e})^{ee''}$ therefore $\{x,v_{l+2}, s_i(v_{l+1})\}$ is also a $2$-simplex and we are done.
	\end{itemize}
 \end{itemize}
\medskip

\textbf{Claim 2.} \textit{The map $\Phi_*\co \pi_1(\Cay(\SAut(\F_n,B),\mathcal{S}),1)\to \pi_1(\lk(B,\PB(\F_n)),v_{l+1})$ is surjective.}

Let $u_0 - u_1 - \ldots - u_k$ be an edge-path loop in $\lk(B,\PB(\F_n))$, with $u_0=u_k=v_{l+1}$. We will show that it is in the image of $\Phi_*$. For $0\leq i <k$, we have that $B\cup \{u_i, u_{i+1}\}$ is a partial basis of $\F_n$.

Next we inductively define elements $\phi_1,\ldots,\phi_k\in \SAut(\F_n,B\cup \{v_{l+1}\})$. Since $B\cup \{v_{l+1},u_1\}$ is a partial basis and $n-l\geq 3$, there is $\phi_1\in \SAut(\F_n,B\cup \{v_{l+1}\})$ such that $\phi_1(v_{l+2})=u_1$. Hence $u_1= \phi_1 w_{v_{l+2},v_{l+1}}(v_{l+1})$.

Now suppose we have defined $\phi_i$ so that $u_i=(\phi_1 w_{v_{l+2},v_{l+1}}) \cdots (\phi_i w_{v_{l+2},v_{l+1}})(v_{l+1})$. Since $B\cup \{u_i,u_{i+1}\}$ is a partial basis, applying the inverse of $(\phi_1 w_{v_{l+2},v_{l+1}}) \cdots (\phi_i w_{v_{l+2},v_{l+1}})$ we see that 
$$B\cup \{v_{l+1},  (\phi_i w_{v_{l+2},v_{l+1}})^{-1} \cdots (\phi_1 w_{v_{l+2},v_{l+1}})^{-1}(u_{i+1})\}$$
is a partial basis and hence there is $\phi_{i+1}\in \SAut(\F_n,B\cup \{v_{l+1}\})$ such that $$\phi_{i+1}(v_{l+2})=(\phi_i w_{v_{l+2},v_{l+1}})^{-1} \cdots (\phi_1 w_{v_{l+2},v_{l+1}})^{-1}(u_{i+1}).$$ Equivalently, $u_{i+1}=(\phi_1 w_{v_{l+2},v_{l+1}}) \cdots (\phi_{i+1} w_{v_{l+2},v_{l+1}})(v_{l+1})$.
We define $$\phi_{k+1}=((\phi_1 w_{v_{l+2},v_{l+1}}) \cdots (\phi_k w_{v_{l+2},v_{l+1}}))^{-1}.$$
Since $(\phi_1 w_{v_{l+2},v_{l+1}}) \cdots (\phi_k w_{v_{l+2},v_{l+1}})(v_{l+1})=u_k=v_{l+1}$, we have $\phi_{k+1}\in  \SAut(\F_n,B\cup \{v_{l+1}\})$. 

For every $1\leq i \leq k+1$, we can find $s_1^i, \ldots , s_{m_i}^i\in \mathcal{S}^{\pm 1}$ that additionally fix $v_{l+1}$ and such that $$\phi_i = s_1^i \cdots s_{m_i}^i.$$
We have
$(s_1^1 \cdots s_{m_1}^1) w_{v_{l+2},v_{l+1}}\cdots  w_{v_{l+2},v_{l+1}}(s_1^{k+1} \cdots s_{m_{k+1}}^{k+1})=1.$
Therefore there is an edge-path loop in $\Cay(\SAut(\F_n),\mathcal{S})$ whose image by $\Phi_*$ is
$$ v_{l+1} - s^1_1(v_{l+1})- s^1_1s^1_2(v_{l+1}) - \cdots - s^1_1\cdots s^1_{m_1}(v_{l+1}) - s^1_1\cdots s^1_{m_1}w_{v_{l+2},v_{l+1}}(v_{l+1}) - \cdots $$
Since $s^i_j(v_{l+1})=v_{l+1}$, for every $1\leq i\leq k+1$ and $1\leq j\leq m_i$, after deleting repeated vertices this path equals
\begin{align*}
v_{l+1} &- (s_1^1s_2^1\cdots s_{m_1}^1 w_{v_{l+2},v_{l+1}})(v_{l+1}) -(s_1^1s_2^1\cdots s_{m_1}^1 w_{v_{l+2},v_{l+1}})(s_1^2s_2^2\cdots s_{m_2}^2 w_{v_{l+2},v_{l+1}})(v_{l+1}) - \cdots \\
&- (s_1^1s_2^1\cdots s_{m_1}^1 w_{v_{l+2},v_{l+1}})(s_1^ks_2^k\cdots s_{m_k}^k w_{v_{l+2},v_{l+1}})(v_{l+1})- v_{l+1}
\end{align*}
which is precisely $u_0 - u_1 - \cdots - u_k$.
\end{proof}
\end{teorema}

\section{A variant of Quillen's result}\label{SeccionQuillen}

If $K$ is a simplicial complex, $\mathcal{X}(K)$ denotes the face poset of $K$.
If $X$ is a poset $\mathcal{K}(X)$ denotes the order complex of $X$. 
The complex $\mathcal{K}(\mathcal{X}(K))$ is the barycentric subdivision $K'$ of $K$.
Throughout the paper we consider homology with integer coefficients and $\wt{C}_\cc(K)$ is the augmented simplicial chain complex.
Let $\lambda:\wt{C}_\cc(K)\to \wt{C}_\cc(K')$ be the subdivision operator $\alpha\mapsto \alpha'$.
If $X$ is a poset we write $\wt{H}_\cc(X)$ for the homology $\wt{H}_\cc(\mathcal{K}(X))$. We thus have $\wt{H}_\cc(X)=\wt{H}_\cc(X^\op)$.
Recall that if $X$ is a poset and $x\in X$ the \textit{height} of $x$ denoted $h(x)$ is the dimension of $\mathcal{K}(X_{\leq x})$.
If $K$ is a simplicial complex we can identify  $\mathcal{X}(\lk(\sigma, K))=\mathcal{X}(K)_{>\sigma}$ by the map $\tau\mapsto \sigma \cup \tau$.
If $K_1,K_2$ are simplicial complexes we have $\wt{C}_\cc(K_1*K_2)=\wt{C}_\cc(K_1)*\wt{C}_\cc(K_2)$ (here $*$ denotes the join of chain complexes, defined as the suspension of the tensor product).
Recall that the join of two posets $X_1,X_2$ is the disjoint union of $X_1$ and $X_2$ keeping the given ordering within $X_1$ and $X_2$ and setting $x_1\leq x_2$ for every $x_1\in X_1$ and $x_2\in X_2$ \cite[Definition 2.7.1]{Barmak2}.
We have $\mathcal{K}(X_1*X_2)=\mathcal{K}(X_1)*\mathcal{K}(X_2)$.
If $X$ is a poset and $x\in X$, then $\lk(x,X)=X_{<x}*X_{>x}$ is the subposet of $X$ consisting of elements that can be compared with $x$. We have $\lk(x,\mathcal{K}(X))=\mathcal{K}(\lk(x,X))$.

\begin{definicion}
Let $f\co X  \to Y$ be an order preserving map. The \textit{non-Hausdorff mapping cylinder $M(f)$} is the poset given by the following order on the disjoint union of $X$ and $Y$. We keep the given ordering within $X$ and $Y$ and for $x\in X$, $y\in Y$ we set $x \leq y$ in $M(f)$ if $f(x) \leq y$ in $Y$.
\end{definicion}
\begin{center}
 \begin{tikzcd}
  X\arrow[hook]{r}[]{j}\arrow[]{dr}[swap]{f}& M(f) \\
  &Y\arrow[hook]{u}[swap]{i}
 \end{tikzcd}
\end{center}

If $j\co X\to M(f)$, $i\co Y\to M(f)$ are the inclusions, then $\mathcal{K}(i)$ is a homotopy equivalence. Since $j\leq if$ we also have $\mathcal{K}(j)\simeq \mathcal{K}(if)$. For more details on this construction see {\cite[2.8]{Barmak2}}.

\begin{definicion}
A simplicial complex $K$ is said to be \textit{$n$-spherical} if $\dim(K)=n$ and $K$ is $(n-1)$-connected.
We say that $K$ is \textit{homologically $n$-spherical} if $\dim(K)=n$ and $\wt{H}_i(K)=0$ for every $i<n$.
Recall that $K$ is \textit{Cohen-Macaulay} if $K$ is $n$-spherical and the link $\lk(\sigma,K)$ is $(n-\dim(\sigma)-1)$-spherical for every simplex $\sigma\in K$.
A poset $X$ is \textit{(homologically) $n$-spherical} if $\mathcal{K}(X)$ is (homologically) $n$-spherical.
\end{definicion}

Recall that if $f\co X\to Y$ is a map of posets, the \textit{fiber of $f$ under $y$} is the subposet $f/y=\{x \tq f(x)\leq y\}\subseteq X$.

\begin{definicion}
An order preserving map $f\co X\to Y$ is \textit{(homologically) $n$-spherical}, if $Y_{>y}$ is (homologically) $(n -h(y)-1)$-spherical and $f/y$ is (homologically) $h(y)$-spherical for all $y\in Y$.
\end{definicion}

\begin{proposicion}\label{cotaobvia}
 Let $f\co X\to Y$ be homologically $n$-spherical. Then for every $x\in X$ we have $h(f(x))\geq h(x)$.
 \begin{proof}
Let $y=f(x)$. Since $x\in f/y$ and $f/y$ is homologically $h(y)$-spherical we have $h(x)\leq \dim(f/y)=h(y)$.
 \end{proof}
\end{proposicion}

\begin{proposicion}
A homologically $n$-spherical map $f\co X\to Y$ is surjective.
 \begin{proof}
 Let $y\in Y$ and let $r=h(y)$. Since $f/y$ is homologically $r$-spherical, $\dim(f/y)=r$. So there is $x\in f/y$ with $h(x)=r$. Let $\wt{y}=f(x)$. We obviously have $\wt{y}\leq y$. By \cref{cotaobvia} we have $h(\wt{y})\geq h(x)=r$. Therefore we have $\wt{y}=y$.
 \end{proof}
\end{proposicion}

From the definition of spherical map we also have the following:
\begin{proposicion}
 If $f\co X\to Y$ is homologically $n$-spherical then $\dim(X)=\dim(Y)=n$.
\end{proposicion}

The first part of the following result is due to Quillen \cite[Theorem 9.1]{Q}. To prove the second part we build on the proof of the first part given by Piterman \cite[Teorema 2.1.28]{P}. The idea of considering the non-Hausdorff mapping cylinder of $f\co X\to Y$ and removing the points of $Y$ from bottom to top is originally due to Barmak and Minian \cite{BM}.

\begin{teorema}\label{QuillenExtendido}
Let $f\co X\to Y$ be a homologically $n$-spherical map between posets such that $Y$ is homologically $n$-spherical. Then $X$ is homologically $n$-spherical, $f_*\co \wt{H}_n(X)\to \wt{H}_n(Y)$ is an epimorphism and
$$\wt{H}_n(X)\simeq \wt{H}_n(Y)\bigoplus_{y\in Y} \wt{H}_{h(y)}(f/y) \otimes \wt{H}_{n-h(y)-1}(Y_{>y}).$$

Moreover suppose that $X=\mathcal{X}(K)$ for certain simplicial complex $K$ and

\noindent\CondB If $f(\sigma_1 )\leq f(\sigma_2)$ then $\lk(\sigma_2,K)\subseteq \lk(\sigma_1,K)$.

\noindent\CondC If $f(\sigma_1)\leq f(\sigma_2)$ and $f(\tau_1)\leq f(\tau_2)$ then $f(\sigma_1\cup \tau_1)\leq f(\sigma_2\cup \tau_2)$, whenever $\sigma_1\cup\tau_1,\sigma_2\cup \tau_2\in K$.

\noindent\CondA For every $y\in Y$ and every $\sigma \in f^{-1}(y)$, the map $f_*\co \wt{H}_{n-h(y)-1}(X_{>\sigma})\to \wt{H}_{n-h(y)-1}(Y_{>y})$ is an epimorphism.

Then we can produce a basis of $\wt{H}_n(K)$ as follows. Since $f_*$ is an epimorphism, we can take $\{\gamma_i\}_{i\in I}\subseteq \wt{H}_n(K)$ such that $\{f_*( \gamma_i')\}_{i\in I}$ is a basis of $\wt{H}_n(Y)$.
In addition, for every $y\in Y$ we choose $x\in f^{-1}(y)$ and we consider the subcomplexes $K_y=\{ \sigma \tq f(\sigma)\leq y\}$ and $K^y=\lk(x,K)$. By \CondB, $K^y$ does not depend on the choice of $x$. Also by \CondB, $K_y*K^y$ is a subcomplex of $K$. Let $\wt{f}\co \mathcal{X}(K^y)\to Y_{>y}$ be defined by $\wt{f}(\tau)=f(x\cup \tau)$. By \CondC, $\wt{f}$ does not depend on the choice of $x$. 
We take a basis $\{\alpha_i \}_{i\in I_{y} }$  of $\wt{H}_{h(y)}(K_y)$ and using \CondA we take $\{\beta_j\}_{j\in J_y}\subseteq \wt{H}_{n-h(y)-1}(K^y)$ such that $\{\wt{f}_*(\beta_j')\}_{j\in J_y}$ is a basis of $\wt{H}_{n-h(y)-1}(Y_{>y})$.
Then $$\{\gamma_i\tq i\in I\}\cup \{ \alpha_i*\beta_j\tq y\in Y, i\in I_y, j\in J_y\}$$ is a basis of $\wt{H}_n(K)$.

\begin{proof}
Let $M=M(f)$ be the non-Hausdorff mapping cilinder of $f$ and let $j\co X\to M$, $i\co Y\to M$ be the inclusions.  We have $j_*=i_*f_*$.
Since $f$ is $n$-spherical we have $\dim(M)=n+1$.

Let $Y_r=\{ y\in Y \tq h(y)\geq r\}$. For each $r$ we consider the subspace $M_r= X\cup Y_r$ of $M$. We have $M_{n+1}=X$ and $M_{0}=M$. Let
$$L_r= \coprod_{h(y)=r} \lk( y, M_r )= \coprod_{h(y)=r} f/y * Y_{>y}.$$

For each $r$ we consider the Mayer-Vietoris sequence for the open covering  $\{U,V\}$ of $\mathcal{K}(M_{r-1})$ given by
\begin{align*}
 U &= \mathcal{K}(M_{r-1})-\{y\in Y\tq h(y)=r-1\}  \\
 V &= \bigcup_{h(y)=r-1} \st(y, \mathcal{K}(M_{r-1})) \\
\end{align*}
where $\st(v,K)$ denotes the open star of $v$ in $K$.
We have homotopy equivalences $U\simeq \mathcal{K}(M_r)$ and $U\cap V\simeq \mathcal{K}(L_r)$. Since $f$ is a homologically $n$-spherical map $\lk( y, M_{r-1})$ is homologically $n$-spherical, so the homology of $L_r$ is concentrated in degrees $0$ and $n$. The tail of the sequence is $0\to \wt{H}_{n+1}(M_r)\to \wt{H}_{n+1}(M_{r-1})$ and since $\wt{H}_{n+1}(M_0)=\wt{H}_{n+1}(Y)=0$ we have $\wt{H}_{n+1}(M_r)=0$ for every $r$. We also have isomorphisms $\wt{H}_i(M_r)\to \wt{H}_i(M_{r-1})$ if $0\leq i \leq n-1$ (since $L_r$ may not be connected, we have to take some care when $i=0,1$). From this we conclude that $X$ is homologically $n$-spherical and we also have short exact sequences

\begin{align*}
0\to \wt{H}_n( L_{n} ) \xrightarrow{i_{n+1}} &\wt{H}_n( X) \xrightarrow{p_{n+1}} \wt{H}_n( M_{n})\to 0\\
&\quad \cdots\\
0\to \wt{H}_n( L_{r-1} ) \xrightarrow{i_r} &\wt{H}_n( M_r) \xrightarrow{p_r} \wt{H}_n( M_{r-1})\to 0\\
&\quad \cdots \\
0\to \wt{H}_n(L_0)  \xrightarrow{i_1}& \wt{H}_n( M_1) \xrightarrow{p_1} \wt{H}_n( M )\to 0.
\end{align*}
Here the map $i_r$ is the map induced by the map $L_{r-1}\to M_{r}$ given by the coproduct of the inclusions $\lk( y, M_{r-1})\to M_{r}$ and the map $p_r$ is induced by the inclusion $M_{r}\to M_{r-1}$. 
By induction on $r$, it follows that these sequences are split and that $\wt{H}_n(M_r)$ is free for every $r$. We have
$$\wt{H}_n(L_r)=\bigoplus_{h(y)=r} \wt{H}_{r}(f/y) \otimes \wt{H}_{n-r-1}( Y_{>y})$$
and therefore using the isomorphism $i_*:\wt{H}_n(Y)\to \wt{H}_n(M)$ we obtain
$$\wt{H}_n(X)=\wt{H}_n(Y)\bigoplus_{y\in Y} \wt{H}_{h(y)}(f/y) \otimes \wt{H}_{n-h(y)-1}(Y_{>y}).$$

Now $\wt{H}_n(j)=p_1\cdots p_n$ is an epimorphism so $f_*\co \wt{H}_n(X)\to \wt{H}_n(Y)$ is also an epimorphism.
We will need the following claim which is proved at the end of the proof.\medskip

\textbf{Claim.} Let $y\in Y$, $r=h(y)$. Then for every $\alpha\in Z_r(K_y)$, $\beta\in Z_{n-r-1}(K^y)$  we have $[ (\alpha*\beta)'] = [\alpha' * \wt{f}_*(\beta')]$ in $\wt{H}_n(M_{r+1})$.\medskip

Let $j_r\co  X\to M_r$ be the inclusion. We have ${j_r}_*=p_{r+1}\circ\ldots\circ p_{n+1}$.
Now by induction on $r$ we prove that for $0\leq r\leq n+1$
$$\{{j_r}_*(\gamma_i')\tq i\in I\}\cup \{ {j_r}_*( (\alpha_i*\beta_j)' ) \tq y\in Y, i\in I_y, j\in J_y, h(y)< r \}$$
is a basis of $\wt{H}_n(M_r)$. Since $j_0=j$ and $j_*=i_*f_*$ it holds when $r=0$.
Now, assuming it holds for $r$, we prove it also holds for $r+1$. By the split exact sequence obtained above, it suffices to check that
$$\{ j_{r+1}( (\alpha_i*\beta_j)' ) \tq i\in I_y, j\in J_y\}$$
is a basis of $\wt{H}_n(\lk(y,M_{r}))$ for every $y\in Y$ of height $r$.
Now in $\wt{H}_{n}(M_{r+1})$ we have
$$j_{r+1}( (\alpha_i*\beta_j)' ) = (\alpha_i*\beta_j)' = \alpha_i' * \wt{f}_*(\beta_j')$$
and the induction is complete, for $\{ \alpha_i'* \wt{f}_*(\beta_j')\}_{i\in I_y, j\in J_y}$ is a basis of $\wt{H}_n( \lk(y,M_r))$. 
We have $j_{n+1}=1_X$ and taking $r=n+1$ we get the desired basis of $\wt{H}_n(K)$.

\begin{proof}[Proof (claim)]
We consider chain maps $\phi_1,\phi_2\co  \wt{C}_\cc( K_y*K^y )\to \wt{C}_\cc(\mathcal{K}(M_{r+1}))$ defined by
 \begin{align*}
  \phi_1\co \wt{C}_\cc( K_y*K^y )&\hookrightarrow \wt{C}_\cc( K)\xrightarrow{\lambda} \wt{C}_\cc(\mathcal{K}(X)) \hookrightarrow  \wt{C}_\cc(\mathcal{K}(M_{r+1}))
 \end{align*}
and
\begin{align*}
 \phi_2\co & \wt{C}_\cc( K_y*K^y )=\wt{C}_\cc(K_y)*\wt{C}_\cc(K^y)\xrightarrow{\lambda*\lambda} \wt{C}_\cc(\mathcal{K}(f/y))*\wt{C}_\cc( \mathcal{K}(\mathcal{X}(K^y)) )\xrightarrow{{1} * \wt{f}_* } \\
 & \quad \wt{C}_\cc(\mathcal{K}(f/y))*\wt{C}_\cc(\mathcal{K}(Y_{>y})) =\wt{C}_\cc( \mathcal{K}(f/y* Y_{>y}) )\hookrightarrow \wt{C}_\cc( \mathcal{K}(M_{r+1})).
\end{align*}
Note that $\phi_1(\alpha*\beta)=(\alpha*\beta)'$ and $\phi_2(\alpha*\beta)= \alpha' * \wt{f}_*(\beta')$. We define an acyclic carrier $\Phi\co K_y*K^y\to \mathcal{K}(M_{r+1})$. If $\sigma\cup \tau$ is a simplex in $K_y*K^y$, with $\sigma\in K_y$ and $\tau\in K^y$, we define $$\Phi(\sigma\cup \tau)=\begin{cases} \mathcal{K}\left({M_{r+1}}_{\leq \wt{f}(\tau)}\right) &\text{ if $\tau\neq \emptyset$. }\\
 \mathcal{K}\left({M_{r+1}}_{\leq \sigma}\right) &\text{if $\tau=\emptyset$.}
 \end{cases}$$
 If $\sigma_1\cup \tau_1\subseteq \sigma_2\cup \tau_2$ are simplices of $K_y*K^y$ where $\sigma_i\in K_y$ and $\tau_i\in K^y$ are possibly empty, we have $\sigma_1\subseteq \sigma_2$ and $\tau_1\subseteq \tau_2$. In $M$ we have $\sigma_1\leq \sigma_2\leq y\leq \wt{f}(\tau_1)\leq \wt{f}(\tau_2)$ so in any case $\Phi(\sigma_1\cup\tau_1)\subseteq\Phi(\sigma_2\cup \tau_2)$. So $\Phi$ is a carrier. It is obviously acyclic.
 
 Now we prove that $\phi_1$ and $\phi_2$ are carried by $\Phi$. To show that $\phi_1$ is carried by $\Phi$ we need to show that $\phi_1(\sigma\cup \tau)=(\sigma\cup \tau)'$ is supported on $\Phi(\sigma\cup \tau)$. If $\tau$ is empty it is clear. If $\tau$ is nonempty, we consider $x\in f^{-1}(y)$.
 In $M$, by \CondC we have $\sigma\cup \tau \leq f(\sigma\cup\tau) \leq f(x\cup\tau)=\wt{f}(\tau)$. Therefore $(\sigma\cup \tau)'$ is supported on $\Phi(\sigma\cup\tau)=\mathcal{K}\left({M_{r+1}}_{\leq \wt{f}(\tau)}\right)$. It is easy to see that $\phi_2$ is also carried by $\Phi$.
 
 Finally by the Acyclic Carrier Theorem \cite[Theorem 13.3]{Munkres} we have $$ [(\alpha*\beta)'] = [\phi_1(\alpha*\beta)]=[\phi_2(\alpha*\beta)]= [ \alpha' * \wt{f}_*(\beta') ]$$ and we are done.
 \end{proof}
\end{proof}
\end{teorema}
\begin{observacion}
 We can consider $\varphi\co X\to M_{r+1}$ given by $$\varphi(x)=\begin{cases}x & \text{if } h(x)<r+1 \\ f(x) & \text{if } h(x)\geq r+1\end{cases}.$$ Then $j_{r+1}\leq \varphi$. Therefore ${j_{r+1}}_*\simeq \mathcal{K}(\varphi)$ and $j_*=\varphi_*$.
In the previous proof we actually have $ \varphi_*( (\alpha*\beta)') = \alpha' * \wt{f}_*(\beta')$ in $Z_n(M_{r+1})$.
\end{observacion}

\section[PB(Fn) is Cohen-Macaulay]{$\PB(\F_n)$ is Cohen-Macaulay}
\label{SeccionPrincipal}

To prove that $\PB(\F_n)$ is Cohen-Macaulay we need to consider other related spaces. The \textit{free factor poset} $\FC(F)$ of a free group $F$ is the poset of proper free factors of $F$ ordered by inclusion. This poset was studied by Hatcher and Vogtmann \cite{HV}. If $H$ is a free factor of $\F_n$ and $B$ is a basis of $H$ then $B$ is a partial basis of $\F_n$. If $B$ is a partial basis of $\F_n$ then $H=\langle B \rangle$ is a free factor of $\F_n$. There is an order preserving map
\begin{align*}
g\co \mathcal{X}\left(\PB(\F_n)^{(n-2)}\right) &\to \FC(\F_n)\\
\sigma &\mapsto \langle \sigma \rangle
\end{align*}
and if $B_0$ is a partial basis we have the restriction $g\co \mathcal{X}\left(\PB(\F_n)^{(n-2)}\right)_{>B_0} \to \FC(\F_n)_{>\langle B_0\rangle}$.

\begin{proposicion}[{\cite[p. 117]{MKS}}]\label{propiedades}
Suppose $H$ is a free factor of $\F_n$ and $K\leq H$. Then $K$ is a free factor of $H$ if and only if $K$ is a free factor of $\F_n$.
\end{proposicion}

\begin{teorema}[Hatcher-Vogtmann, {\cite[\S 4]{HV}}]\label{ResultadoHV}
If $H\leq \F_n$ is a free factor, $\FC(\F_n)_{>H}$ is $(n-\rk{H}-2)$-spherical.
\end{teorema}

We will consider the following simplicial complex $Y$ with vertices the free factors of $\F_n$ that have rank $n-1$. A simplex of $Y$ is a set of free factors $\{ H_1,\ldots, H_k\}$ such that there is a basis $\{w_1,\ldots,w_n\}$ of $\F_n$ such that for $1\leq i\leq k$ we have $H_i=\langle w_1,\ldots, w_{i-1},w_{i+1},\ldots,w_n\rangle$.  If $H\leq \F_n$ is a free factor, we consider the full subcomplex $Y_H$ of $Y$ spanned by the vertices which are free factors containing $H$. There is another equivalent definition for $Y$ and $Y_H$ in terms of sphere systems, see \cite[Remark after Corollary 3.4]{HV}.

\begin{teorema}[Hatcher-Vogtmann,{\cite[Theorem 2.4]{HV}}]
 Let $H$ be a free factor of $\F_n$. Then $Y_H$ is $(n-\rk{H}-1)$-spherical.
\end{teorema}
There is a spherical map $f\co \mathcal{X}(Y_H^{(n-\rk{H}-2)})\to (\FC(\F_n)_{> H})^{\op}$ that maps a simplex $\sigma=\{H_1,\ldots, H_k\}$ to $H_1\cap \cdots \cap H_k$. Hatcher and Vogtmann used the map $f$ to prove \cref{ResultadoHV}. We also consider the map $\wt{g}\co \mathcal{X}\left(\lk\left( B, \PB(\F_n)^{(n-2)}\right)\right) \to \FC(\F_n)_{>\langle B \rangle}$ given by  $\sigma \mapsto \langle B\cup \sigma \rangle$ which can be identified with $g\co \mathcal{X}\left(\PB(\F_n)^{(n-2)}\right)_{>B} \to \FC(\F_n)_{>\langle B \rangle} $. The following technical lemma will be needed later.

\begin{lema}\label{LemaNuevo}
Let $B$ be a partial basis of  $\F_n$, $|B|=l$. Let $\overline{\gamma} \in \wt{H}_{n-l-2}(\FC(\F_n)_{>\langle B\rangle})$. There exists $\gamma\in B_{n-l-2}(\lk(B,\PB(\F_n)))$ such that $\wt{g}_*( \gamma' )=\overline{\gamma}$.
\begin{proof}
We define a map $\phi\co C_{n-l-1}(Y_{\langle B\rangle})\to C_{n-l-1}(\lk(B,\PB(\F_n)))$ as follows. For each $(n-l-1)$-simplex $\sigma=\{H_{l+1},\ldots,H_{n}\}$ of $Y_{\langle B\rangle }$ we choose a basis $\{w_1,\ldots, w_n\}$ of $\F_n$ such that $$H_i=\langle w_1,\ldots, w_{i-1},w_{i+1},\ldots, w_n\rangle$$ for $l+1\leq i \leq n$. Then by \cref{propiedades} we have $$\langle B\rangle = \bigcap_{i=l+1}^n H_i = \langle w_{1}, \ldots w_l\rangle$$ so $B\cup \{w_{l+1},\ldots, w_{n}\}$ is a basis of $\F_n$. Therefore $\wt{\sigma} = \{w_{l+1},\ldots, w_{n}\}$ is an $(n-l-1)$-simplex of $\lk(B, \PB(\F_n))$. Then we define the map $\phi$ on $\sigma$ by $\phi(\sigma)=\wt{\sigma}$.

Now since $f\co \mathcal{X}(Y_{\langle B\rangle}^{(n-l-2)}) \to  (\FC(\F_n)_{>\langle B\rangle})^{\op}$ is $(n-l-2)$-spherical \cite[\S 4]{HV}, by \cref{QuillenExtendido} we have an epimorphism $f_*\co  \wt{H}_{n-l-2} (\mathcal{X}({Y_{\langle B\rangle}^{(n-l-2)}}))\to  \wt{H}_{n-l-2}(\FC(\F_n)_{>\langle B\rangle})$ and since $Y_{\langle B\rangle}$ is $(n-l-2)$-connected, there is $c\in C_{n-l-1}(Y_{\langle B\rangle})$ such that $f_*( d(c)' )=\overline{\gamma}$ . We define $\gamma=d\phi(c)$. We immediately have $\gamma\in B_{n-l-2}(\lk(B,\PB(\F_n)))$. It is easy to verify that $\wt{g}_*( d\phi(\sigma)')=f_*(d\sigma')$ and from this it follows that $\wt{g}_*(\gamma')=\wt{g}_*(d\phi(c)')=f_*(dc')=\overline{\gamma}$.
\end{proof}
\end{lema}
 
\begin{teorema}
The complex $\lk(B_0, \PB(\F_n))$ is $(n-|B_0|-1)$-spherical for any partial basis $B_0$ of $\F_n$.
 \begin{proof} 
We proceed by induction on $k=n-|B_0|$. If  $k\leq 3$ it follows from \cref{linksimplementeconexo}. Now if $k\geq 4$ we want to apply \cref{QuillenExtendido} to the map
$g\co \mathcal{X}\left(\PB(\F_n)^{(n-2)}\right)_{>B_0} \to \FC(\F_n)_{>\langle B_0\rangle}$.

By \cref{ResultadoHV}, $\FC(\F_n)_{>\langle B_0\rangle}$ is $(n-|B_0|-2)$-spherical.
In addition $g$ is spherical, since $\FC(\F_n)_{>H}$ is $(n-\rk{H}-2)$-spherical if $H\in \FC(\F_n)_{>\langle B_0\rangle}$ and by the induction hypothesis $g/H=\mathcal{X}\left(\PB(H)\right)_{>B_0}=\mathcal{X}(\lk(B_0,\PB(H)))$ is $(\rk{H}-|B_0|-1)$-spherical.
Then by \cref{QuillenExtendido}, $\mathcal{X}\left(\PB(\F_n)^{(n-2)}\right)_{>B_0}$ is homologically $(n-|B_0|-2)$-spherical. 

We identify $\mathcal{X}\left(\PB(\F_n)^{(n-2)}\right)_{>B_0} = \mathcal{X}( \lk( B_0, \PB(\F_n)^{(n-2)}) )$.
Now we check the hypotheses \CondB, \CondC and \CondA of \cref{QuillenExtendido}.
If $\wt{g}(B_1)\subseteq \wt{g}(B_2)$ it is easy to see that $\lk(B_2, \lk( B_0, \PB(\F_n)^{(n-2)}))\subseteq \lk(B_1,\lk( B_0, \PB(\F_n)^{(n-2)}))$ so \CondB holds.
Obviously \CondC holds.
And by the induction hypothesis \CondA holds.
Thus, the second part of  \cref{QuillenExtendido} gives a basis of $\wt{H}_{n-|B_0|-2}(\lk( B_0, \PB(\F_n)^{(n-2)}))$. By \cref{LemaNuevo} we can choose the $\gamma_i$ to be borders.
We need to prove that the remaining elements of this basis are trivial in $\wt{H}_{n-|B_0|-2}(\PB(\F_n)_{>B_0})$. We only have to show that for all $H\in \FC(\F_n)_{>\langle B_0 \rangle}$, $i\in I_H$, $j\in J_H$
$$\alpha_i*\beta_j\in B_{n-|B_0|-2}(\lk(B_0, \PB(\F_n)))$$
We take a basis $B$ of $H$. By the induction hypothesis we have $\wt{H}_{n-|B|-2}(\lk( B, \PB(\F_n)))=0$. So there is $\omega\in C_{n-|B|-1}(\lk(B,\PB(\F_n)))$ such that $d(\omega)=(-1)^{|\alpha_i|}\, \beta_j$. Therefore
$$d( \alpha_i* \omega  )= d(\alpha_i) * \omega + (-1)^{|\alpha_i|}\, \alpha_i*d(\omega) = \alpha_i*\beta_j. $$

Therefore  $\lk(B_0, \PB(\F_n))$ is homologically $(n-|B_0|-1)$-spherical and by \cref{linksimplementeconexo} it is $(n-|B_0|-1)$-spherical.
\end{proof}
\end{teorema}

\begin{observacion}
\cref{QuillenExtendido} also holds without the word \textit{homologically} (see \cite[Theorem 9.1]{Q}). Thus, we may easily modify the previous proof so that \cref{linksimplementeconexo} is only used as the base case $k\leq 3$.
\end{observacion}

\begin{corolario}\label{PBFnCohenMacaulay}
The complex $\PB(\F_n)$ is Cohen-Macaulay of dimension $n-1$.
\begin{proof}
It follows immediately from the previous result.
\end{proof}
\end{corolario}

\end{document}